\documentclass[11pt]{article}
\usepackage{amsmath,amssymb,amscd,enumerate,eucal,bm}
\numberwithin{equation}{section}
\newcommand{\R}{\mathbb{R}}

\begin{document}
\setlength{\baselineskip}{15pt}
\title{Norm Constants in cases of the Caffarelli-Kohn-Nirenberg Inequality}
\author{Akshay L. Chanillo$^{(1)}$; Sagun Chanillo$^{(2)}$ and Ali Maalaoui$^{(3)}$}

\addtocounter{footnote}{1}
\footnotetext{20 Dorland Farm Ct., Skillman, NJ 08558, USA. E-mail address:
{\tt{achanillo3@gmail.com}}}

\addtocounter{footnote}{1}
\footnotetext{Deptt. of Math., Rutgers Univ. 110 Frelinghuysen Rd., Piscataway, NJ 08854 USA. E-mail address:
{\tt{chanillo@math.rutgers.edu}}}

\addtocounter{footnote}{1}
\footnotetext{Department of mathematics and natural sciences, American University of Ras Al Khaimah, PO Box 10021, Ras Al Khaimah, UAE. E-mail address:
{\tt{ali.maalaoui@aurak.ac.ae}}}

\date{}
\maketitle
\begin{abstract}
Based on elementary Linear algebra, we
provide radically simplified proofs using quasi-conformal changes of
variables to obtain sharp constants and optimizers in cases of the
Caffarelli-Kohn-Nirenberg inequality. Some of our results were
obtained earlier by Lam and Lu.
{\footnotesize \baselineskip 4mm }
\end{abstract}


\medskip

{\footnotesize\textit{}}

\newtheorem{theorem}{Theorem}[section]
\newtheorem{lemma}{Lemma}[section]
\newtheorem{corollary}{Corollary}[section]
\newtheorem{definition}{Definition}[section]
\newtheorem{proposition}{Proposition}[section]
\newtheorem{remark}{Remark}[section]
\newtheorem{claim}{Claim}[section]

\section{Introduction}

The aim of this elementary note is to see if we can obtain sharp constants and also the minimizer for special cases of the
Cafarelli-Kohn-Nirenberg inequality via change of variables employing a suitable Quasi-conformal map. We start by recalling the main inequality that we are interested in proved first in \cite{CKN}:
\begin{equation}\label{CKN}
\bigg(\int_{\R^n} |f(x)|^r\ |x|^{\alpha n} dx\bigg)^{1/r}\leq C \bigg(\int_{\R^n}|f(x)|^s\ |x|^{\alpha n}\ dx\bigg)^{{{1-t}\over{s}}}\bigg(\int_{R^n}|\nabla f(x)|^p\ |x|^{\alpha(n-p)}\
   dx\bigg)^{{{t}\over{p}}},
\end{equation}
where $1\leq p<n$ and $\frac{1}{r}=\frac{1-t}{s}+\frac{t(n-p)}{np}$ with $0\leq t\leq 1$ and $a>-1$.

In our investigation of the best constant $C$ for which inequality $(\ref{CKN})$ holds, we will use the constant $M(s,r,p)$, for $s$, $r$ and $p$ as above, appearing in the inequality
\begin{equation}\label{dd}
\|f\|_{r}\leq M(s,r,p)\|f\|^{1-t}_{s}\|\nabla f\|_{p}^{t}.
\end{equation}
The constant and optimizers of this last inequality were investigated by Del Pino and Dolbeault in \cite{DD}.
In fact the proof of $(\ref{dd})$ follows by noticing that
            $$\|f\|_r\leq \|f\|_s^{1-t}\|f\|^t_{{{np}\over{n-p}}}.$$
We then apply the Sobolev inequality to the second term on the right
               $$\|f\|_{{{np}\over{n-p}}}\leq C_p \|\nabla f\|_p,$$
to conclude. In this paper, we show that the optimization problem
for the best constant of $(\ref{CKN})$ exhibit different behaviour
for the two cases $\alpha>0$ and $0>\alpha>-1$. In the case
$0>\alpha>-1$, we compute the best constant and we show that the
optimizer is radial. In the case $\alpha>0$, we also compute the
best constant and we show that in this case there is a break in the
symmetry of the optimizers, since the best constant cannot be
obtained by radial functions anymore. More importantly we establish
that there is no optimizer. We mainly rely on a study of the
eigenfunctions and eigenvalues of the differential of the
quasi-conformal change of variable that we will use. The results of
this paper can be stated as follows
\begin{theorem}
The sharp constant in the CKN inequality $(\ref{CKN})$ for $-1<\alpha<0$ and any function $f$, radial or otherwise,
is given by, $$(1+\alpha)^{{{t}\over{n}}-t}M(s,r,p),$$ where $M(s,r,p)$ is the constant in $(\ref{dd})$. Moreover the optimizer of $(\ref{CKN})$ is then a radial function and can be taken to be $f\circ\phi$, where $f$ can be taken to be the radial optimizers in the cases investigated in [DD].
\end{theorem}
The result in Theorem 1.1 was established earlier by Nguyen Lam and Guozhen Lu \cite{LL} using the quasi-conformal map that we use in our work. Our proof in part is a radically simplified approach. A very comprehensive list of references on this topic is found in \cite{LL} and also \cite{DEL}. In addition to Theorem 1.1 we prove a symmetry breaking phenomena for the case $\alpha>0$. Indeed, we have
\begin{theorem}
The sharp constant for the CKN for $\alpha>0$ is given by
$$(1+\alpha)^{\frac{t}{n}}M(s,r,p).$$
Moreover, there is no optimizer for the inequality and this constant is strictly bigger than the one obtained for radial functions.
\end{theorem}
Our proofs also show the following theorem for radial functions for
all $\alpha>-1$.
\begin{theorem}
The sharp constant for inequality $(\ref{CKN})$ restricted to radial functions for $\alpha>-1$ is given by
$$(1+\alpha)^{{{t}\over{n}}-t}M(s,r,p).$$
\end{theorem}

\section{Proof of the Theorems}
Our aim is now to make quasi-conformal changes of variables in $(\ref{dd})$ with explicit information on the Jacobian and eigenvalues of the quasi-conformal changes of variables. To this end the eigenvalues and Jacobians can be calculated explicitly using a Linear Algebra trick in Lemma (2.23, pg. 14)  in the paper by Chanillo-Torchinsky \cite{CT}. This method was introduced in \cite{CT} to calculate the Hessian that appears in stationary phase calculations, but the Linear Algebra idea goes back to the calculation of what are called permanents and bordered matrices and may be found in classical books on Algebra.\\

If we set $y=\phi(x)$ in $(\ref{dd})$, then
                   $$Dy=D\phi(x)$$
and so if we set $A=D\phi$, $(\ref{dd})$ becomes using the chain rule,

          \begin{align}
 \bigg(\int_{\R^n} |f\circ \phi(x)|^r |J_\phi(x)|\  dx\bigg)^{1/r}\leq& M(s,r,p)\bigg(\int_{\R^n}|f\circ \phi(x)|^s|J_\phi(x)|\
             dx\bigg)^{(1-t)/s}\notag \\
&\times \bigg(\int_{\R^n}|(A^{-1})^\star\nabla(f\circ\phi)|^p|J_\phi(x)|\
             dx\bigg)^{t/p}.\label{eq1}
\end{align}
Here $J_\phi$ is the Jacobian of the map $y=\phi(x)$ that is
$|J_\phi(x)|=|\det D\phi(x)|$ and $B^\star$ denotes the transpose of
the matrix $B$. To obtain the Caffarelli-Kohn-Nirenberg inequalities
in some cases we simply choose the explicit Quasi-conformal map (see
also \cite{CW}):
               \begin{equation}\label{phi}
\phi(x)=x|x|^\alpha, \alpha>-1.
\end{equation}
Thus the goal is to calculate explicitly the eigenvalues of the differential of $(\ref{phi})$ and thus we have full information of the
matrix $A$ above and in particular the Jacobian of $(\ref{phi})$.

\begin{remark}
We also remark that using the work \cite{DD} we can also consider the case of $p=n$ and the Onofri inequality and
Moser-Trudinger type inequalities.
\end{remark}
\begin{lemma}
 Given the map $\phi(x)$ as in $(\ref{phi})$, the differential $A=D\phi(x)$ is unitarily diagonalizable and the
eigenvalues of $A$ are given by
     $$\lambda_1=(1+\alpha)|x|^\alpha,\ \lambda_2=\cdots
     =\lambda_n=|x|^\alpha.$$
Thus as a corollary we obtain that the Jacobian
                  $$|J_\phi(x)|=(1+\alpha)|x|^{\alpha n}.$$
\end{lemma}

{\it Proof:} The proof of Lemma 2.1 involves implementing the elementary proof of Lemma 2.23, in \cite{CT} in this special
situation. Since $D\phi$ is a symmetric matrix, $D\phi$ is unitarily diagonalizable, that is one can write $A=QRQ^t$ where $Q$ is a rotation matrix and $R$ a diagonal matrix. It is enough to compute the eigenvalues for $D\phi$, the Jacobian formula follows by multiplying the eigenvalues. First note,
        $$D\phi(x)=A=\begin{pmatrix}
 |x|^\alpha+\alpha x_1^2|x|^{\alpha-2}& \cdots& \alpha x_1x_j|x|^{\alpha-2}&\cdots &\alpha x_1x_n|x|^{\alpha-2}\cr
        \vdots& &\vdots& &\vdots \cr
        \alpha x_1x_n|x|^{\alpha-2}&\cdots&\alpha x_nx_j|x|^{\alpha-2}&\cdots&
        |x|^\alpha+\alpha x_n^2|x|^{\alpha-2}\cr \end{pmatrix}.$$
Next
        $$A-\lambda I=\begin{pmatrix}
 |x|^\alpha+\alpha x_1^2|x|^{\alpha-2}-\lambda& \cdots& \alpha x_1x_j|x|^{\alpha-2}&\cdots &\alpha x_1x_n|x|^{\alpha-2}\cr
        \vdots& &\vdots& &\vdots\cr
        \alpha x_1x_n|x|^{\alpha-2}&\cdots&\alpha x_nx_j|x|^{\alpha-2}&\cdots&
        |x|^\alpha+\alpha x_n^2|x|^{\alpha-2}-\lambda\cr
\end{pmatrix}.$$

It follows that
           $$\det(A-\lambda I)=|x|^{n(\alpha-2)}\det C,$$
where
           $$C=\begin{pmatrix}
 |x|^2+\alpha x_1^2-\lambda|x|^{2-\alpha} & \cdots&\alpha x_1x_j&\cdots &\alpha x_1x_n\cr
        \vdots& &\vdots& &\vdots\cr
        \alpha x_1x_n&\cdots& \alpha x_nx_j&\cdots&
        |x|^2+\alpha x_n^2-\lambda |x|^{2-\alpha}\cr\end{pmatrix}.$$
To compute the characteristic polynomial of $D\phi$ we simply
compute $\det C$. It is now that we use the trick in \cite{CT}. We simply
add an extra row and column to $C$ such that the new matrix now with
$n+1$ rows and $n+1$ columns has the same determinant as $C$. Thus
we form the matrix $D$, given by

$$D=\begin{pmatrix}
 1&x_1&\cdots& x_j&\cdots &x_n\cr
0&|x|^2+\alpha x_1^2-\lambda|x|^{2-\alpha} & \cdots&    \alpha x_1x_j&\cdots &x_1x_n\cr
\vdots&\vdots& &\vdots&\vdots&\vdots\cr
 0&\alpha x_1x_n&\cdots&\alpha x_nx_j&\cdots&
        |x|^2+\alpha x_n^2-\lambda |x|^{2-\alpha}\cr\end{pmatrix}.$$
Note $\det C=\det D$. Now we perform the elementary row operation in
$D$ that preserves the determinant. We replace row $R_j$, $j\geq 2$
by $R_j-\alpha x_{j-1}R_1$ where $R_1$ is row 1. The new matrix we
get is

$$E=\begin{pmatrix}
1&x_1&\cdots& x_j&\cdots &x_n\cr
-\alpha x_1&|x|^2-\lambda|x|^{2-\alpha} & \cdots&  0&\cdots & 0\cr
\vdots&\vdots & &\vdots & & \vdots \cr
-\alpha x_j&0&\cdots& |x|^{2}-\lambda|x|^{2-\alpha}&\cdots &0\cr
\vdots&\vdots & &\vdots & & \vdots \cr
        -\alpha x_n&0 &\cdots&0&\cdots&
        |x|^2-\lambda |x|^{2-\alpha}\cr\end{pmatrix}.$$
Note $\det D=\det E$ and the matrix we get if we remove the first
row and first column of $E$ is a diagonal matrix with $|x|^2-\lambda
|x|^{2-\alpha}$ on the diagonal. To compute $\det E$ we simply
expand for the determinant using the first row of $E$ and then
expanding by the $j$-th row of the subsequent co-factor matrices we
get for the entry $a_{1,j+1}$. We get,
        $$\det E= (|x|^2-\lambda |x|^{2-\alpha})^n+\alpha(|x|^2-\alpha
        |x|^{2-\alpha})^{n-1}\sum_{j=1}^nx_j^2,$$
which is
$$(|x|^2-\lambda |x|^{2-\alpha})^n+\alpha(|x|^2-\alpha
        |x|^{2-\alpha})^{n-1}|x|^2.$$
The expression above obviously factors as
             \begin{equation}\label{eig}
 (|x|^2-\lambda |x|^{2-\alpha})^{n-1}(
             (1+\alpha)|x|^2-\lambda |x|^{2-\alpha}).
\end{equation}
From $(\ref{eig})$ the conclusion of our Lemma follows because
         $$\det(A-\lambda I)=0=(|x|^2-\lambda |x|^{2-\alpha})^{n-1}(
             (1+\alpha)|x|^2-\lambda |x|^{2-\alpha}).$$
\hfill$\Box$

\medskip
First note that $A=D\phi$ is a symmetric matrix and thus there exist
rotation matrices $Q$ such that,
          $$ A^{-1}=(A^{-1})^\star=Q\hat{D}Q^t,$$
where $\hat{D}$ is diagonal. Using the eigenvalues of $A$ computed
from  Lemma 2.1 above we may write
        $$\hat{D}=|x|^{-\alpha} D,$$
where
 $$D={\rm diag}( (1+\alpha)^{-1},1,1,\cdots, 1).$$

 Thus, by Lemma 2.1,

$$\int_{\R^n}|(A^{-1})^\star\nabla f|^p|J_\phi(x)|\
             dx=(1+\alpha)\int_{\R^n} |QDQ^t(\nabla f)|^p |x|^{\alpha(n-p)}\ dx. $$
We now apply Lemma 2.1 to $(\ref{eq1})$ and we get with $M(s,r,p)$ the constant that occurs in $(\ref{dd})$,
   \begin{align}
\bigg(\int_{\R^n} |f(x)|^r\ |x|^{\alpha n} dx\bigg)^{1/r}&\leq(1+\alpha)^{{{t}\over{n}}} A_\alpha
   M(s,r,p) \bigg(\int_{\R^n}|f(x)|^s\ |x|^{\alpha n}\ dx\bigg)^{{{1-t}\over{s}}}\notag \\
   &\times\bigg(\int_{\R^n}|\nabla f(x)|^p\ |x|^{\alpha(n-p)}\  dx\bigg)^{{{t}\over{p}}},\label{eq2}
   \end{align}
where we define
    $$A_\alpha=\sup_f \bigg[{{\int_{\R^n}|QDQ^t(\nabla
    f)(x)|^p|x|^{\alpha(n-p)}\ dx}\over{\int_{\R^n}|\nabla
    f(x)|^p|x|^{\alpha(n-p)}\ dx}}\bigg]^{{{t}\over{p}}}.$$
The supremum is taken over those functions $f$ where the denominator
in the definition above is finite.

\bigskip
\begin{lemma}
For $\alpha>-1$,
        $$B_\alpha\leq A_\alpha\leq C_\alpha=\begin{cases}1, \alpha\geq 0\cr
                              (1+\alpha)^{-t},-1<\alpha<0\cr\end{cases},$$
   with,
            $$B_\alpha=\bigg[{{\int_{S^{n-1}} \big( [
            {{1}\over{(1+\alpha)^2}}-1]\cos^2\psi(\sigma)+1\big)^{{{p}\over{2}}}\
            d\sigma}\over{\int_{S^{n-1}}
            d\sigma}}\bigg]^{{{t}\over{p}}},$$
   where,
           $$\cos\psi(\sigma)=<\sigma,v>,$$
   where $v$ is a unit eigenvector for the eigenvalue
   $1/(1+\alpha)$ at $\sigma \in S^{n-1}$.
\end{lemma}
On radial functions $f$ we may take $A_\alpha=B_\alpha$ for any $\alpha>-1$ and the ratio defining $A_\alpha$ is identically
$B_\alpha$ for all radial functions without the supremum.

\bigskip

{\it Proof:} We now verify the assertions made about $A_\alpha$. We note that pointwise
          $$|QDQ^t(\nabla f)(x)|=|DQ^t\nabla f(x)|\leq C_\alpha^{{1}\over{t}} |\nabla f(x)|.$$
This establishes $A_\alpha\leq C_\alpha.$

Next we establish the lower bound on $A_\alpha$. Here we assume $f$
is radial. Now notice
      \begin{equation}\label{A}
 |x|^{-\alpha}A=QD^{-1}Q^t
     \end{equation}
the coefficients of $A$ are homogeneous of degree $\alpha$ and thus the coefficients of the left side of $(\ref{A})$ are homogeneous of degree $0$ and since $D$ is a constant matrix, it follows that the
coefficients of $Q, Q^t$ are functions of $\sigma\in S^{n-1}$. We now wish to consider for $f$ radial the expression,
     $${{\int_{\R^n} |QDQ^t(\nabla f)|^p |x|^{\alpha(n-p)}\ dx}\over{\int_{\R^n} |\nabla f|^p |x|^{\alpha(n-p)}\
     dx}}.$$
Using the fact that the coefficients of $Q$ depend only on $\sigma$
we see when $f$ is radial, the expression above when converted to
polar coordinates is equal to:
              \begin{equation}\label{AA}
 {{\int_{S^{n-1}}|Q(\sigma)DQ^t(\sigma)\nabla r|^p\
              d\sigma}\over{\int_{S^{n-1}}d\sigma}}.
\end{equation}
Now let $\{e_i(\sigma)\}_{i=1}^n$ be an orthonormal basis of eigenvectors for $QDQ^t$. Then since $\nabla r=\sigma$ we see that
  $$ |QDQ^t\nabla r|^2=
  {{1}\over{(1+\alpha)^2}}<e_1,\sigma>^2+\sum_{j=2}^n
  <\sigma,e_j>^2$$
The expression above can be re-arranged obviously as
 $$ \big( {{1}\over{(1+\alpha)^2}}-1\big)\cos^2\psi+1.$$
Substituting this expression into $(\ref{AA})$ we readily establish
           $$B_\alpha\leq A_\alpha.$$
Note since we have equality at every step in the computation above,
we also obtain that $A_\alpha=B_\alpha$ when $f$ is radial.

\hfill$\Box$

\begin{lemma}
 For $\alpha>-1$, we have,
          $$B_\alpha=(1+\alpha)^{-t}.$$
\end{lemma}

{\it Proof:}   From the formula for $A=D\phi$, the eigenvalue equation for $A-(1+\alpha)I$ is
     \begin{equation}\label{eig2}
(\sigma_j^2-1)y_j + \sum_{k\not=j}\sigma_k\sigma_j y_k=0,\ j=1,2,\cdots,n
\end{equation}
where the eigenvector is $v=y=(y_1,y_2,\cdots, y_n)$. Now we set $y=\sigma=(\sigma_1,\sigma_2,\cdots ,\sigma_n)$ and we get, the left side of $(\ref{eig2})$ is
        $$\sigma_j(\sigma_j^2-1)+\sigma_j\sum_{k\not
        =j}\sigma_k^2=\sigma_j(\sigma_j^2-1)+\sigma_j(1-\sigma_j^2)=0.$$
Thus $\sigma$ is the unit eigenvector for the eigenvalue
$(1+\alpha)$, which we already know has a 1-dimensional eigenspace.
Thus,
        $$\cos \psi=<\sigma,\sigma>=1,$$
and it follows from the expression for $B_\alpha$ in the statement of Lemma 2.2 that for any $\alpha>-1$,
$$B_\alpha=(1+\alpha)^{-t}.$$

\hfill$\Box$

Notice that this Lemma shows in particular that if we restrict $(\ref{CKN})$ to radial functions, then the sharp constant is $(1+\alpha)^{\frac{t}{n}-t}M(s,r,p),$ as stated in Theorem 1.3.
\bigskip\noindent
\begin{corollary} When $-1<\alpha<0$, then
        $$B_\alpha=A_\alpha=C_\alpha=(1+\alpha)^{-t}.$$
\end{corollary}
\bigskip
{\it Proof:} The proof is obvious combining the conclusions of Lemma 2.2 and Lemma 2.3, which yields
 $B_\alpha=C_\alpha$
when $-1<\alpha<0$.

\hfill$\Box$.

\bigskip
Thus the sharp constant in $(\ref{eq2})$  is established when $-1<\alpha<0$.\\
\section{The case $\alpha>0$}
For the case $\alpha>0$, so far we have the upper and lower bounds
for the sharp constant that is   $(1+\alpha)^{-t}\leq A_\alpha\leq
1,$ in $(\ref{eq2})$. Also if we restrict to radial functions then
$$A_{\alpha}(radial)=(1+\alpha)^{-t}.$$ Now if we check closely the
computations in Lemma 2.3, we have that
 $$|QDQ^t(\nabla f)(x)|^2=
 \big[\big({{1}\over{(1+\alpha)^2}}-1\big)\cos^2\psi+1)\big]|\nabla
 f(x)|^2,$$
where $\cos\psi=<v,w>$, $v$ is a unit vector in the direction of the eigenvector for $(1+\alpha)$ and $w$ the unit vector in the
direction of $\nabla f$. In particular $v$ is radial at all points $x\in \R^n$. But now for $\alpha>0$, $(1+\alpha)^{-2}-1<0$
and so it is advantageous to arrange $\cos\psi=0$ as opposed to $\alpha<0$ when $(1+\alpha)^{-2}-1>0$ and so there
it is advantageous to have $\cos \psi=1$ or functions to be radial. So the idea of proving Theorem 2, is to choose a function gradient having a big angular component that dominates the radial component. If one wants an optimizer for the case $\alpha>0$,  $\nabla f$ needs to be orthogonal to $v$, that is tangent to the sphere at all points. Notice the tangential directions to the sphere are eigenvectors to the eigenvalue $1$ for $A=D\phi$. But in particular in 3D, $\nabla\times \nabla f=0$, no such functions exist, or if $f$ has some smoothness the vector field $\nabla f$ on $S^2$ will be smooth and tangential to $S^2$ which cannot happen by the Hairy ball theorem. In fact, we show the following
\begin{lemma}
If $\alpha>0$, then $A_{\alpha}=1$.
\end{lemma}
{\it Proof:} Indeed, based on the computations above, we have that
the matrix $\tilde{A}=QDQ^{t}$ has two eigenvalues. The first one is
$\frac{1}{(1+\alpha)}<1$, corresponding to the radial direction
$\nabla r$ and the second eigenvalue is $1$ with multiplicity
$(n-1)$ corresponding to the angular directions (tangential to
$S^{n-1}$). We want to estimate the quantity
$$F(f)=\frac{\int_{\mathbb{R}^{n}}|\tilde{A}\nabla f|^{p}|x|^{\alpha(n-p)}dx}{\int_{\mathbb{R}^{n}}|\nabla f|^{p}|x|^{\alpha(n-p)}dx}$$
for some choice of function $f$ knowing that
$A_{\alpha}=\sup_{f}F(f)$.  We use spherical coordinates
$(r,\phi_{1},\cdots,\phi_{n-1})$, and we form
$$f_{k}(r,\phi_{1},\cdots,\phi_{n-1})=h(r)\sin(\phi_{1})\cdots
\sin(\phi_{n-2})\cos(k\phi_{n-1}),$$ where $h:[0,\infty)\to
\mathbb{R}$ is smooth and $h(t)=0$ for $t<1$ and $t>4$.
For the sake of simplicity, we do the computation in $n=3$, the higher dimensional case is similar.\\
 So $f_{k}=h(r)\sin(\phi)\cos(k\theta)$, thus
$$\nabla f_{k}=h'(r)\sin(\phi)\cos(k\theta)\boldsymbol{u}_{r}+\frac{h(r)}{r}\cos(\phi)\cos(k\theta)\boldsymbol{u}_{\phi}-\frac{h(r)}{r}k\sin(k\theta)\boldsymbol{u}_{\theta},$$
where $(\boldsymbol{u}_{r}, \boldsymbol{u}_{\phi}, \boldsymbol{u}_{\theta})$ is the standard orthonormal base defining the spherical coordinate system. Thus $$\tilde{A}\nabla f_{k}=\frac{h'(r)\sin(\phi)\cos(k\theta)}{(1+\alpha)}\boldsymbol{u}_{r}+\frac{h(r)}{r}\cos(\phi)\cos(k\theta)\boldsymbol{u}_{\phi}-\frac{h(r)}{r}k\sin(k\theta)\boldsymbol{u}_{\theta}$$
We compute then
\begin{align}
|\tilde{A}\nabla f_{k}|^{p}&=\left[\frac{h'(r)^{2}\sin^{2}(\phi)\cos^{2}(k\theta)}{(1+\alpha)^{2}}+\left(\frac{h(r)}{r}\right)^{2}\cos^{2}(\phi)\cos^{2}(k\theta)+\left(\frac{h(r)}{r}\right)^{2}k^{2}\sin^{2}(k\theta)\right]^{\frac{p}{2}}\notag\\
&=\left[\cos^{2}(k\theta)(\frac{h'(r)^{2}\sin^{2}(\phi)}{(1+\alpha)^{2}}+\left(\frac{h(r)}{r}\right)^{2}\cos^{2}(\phi))+k^{2}\left(\frac{h(r)}{r}\right)^{2}\sin^{2}(k\theta)\right]^{\frac{p}{2}}\notag
\end{align}
Hence,
\begin{align}
\int_{\mathbb{R}^{3}}|\tilde{A}\nabla f_{k}|^{p}|x|^{\alpha(n-p)}\ dx&=\int_{1}^{4}\int_{0}^{\frac{\pi}{2}}\int_{0}^{2\pi}\Bigg[\cos^{2}(k\theta)\left(\frac{h'(r)^{2}\sin^{2}(\phi)}{(1+\alpha)^{2}}+\left(\frac{h(r)}{r}\right)^{2}\cos^{2}(\phi)\right)\notag\\
&\qquad \qquad \qquad+k^{2}\left(\frac{h(r)}{r}\right)^{2}\sin^{2}(k\theta)\Bigg]^{\frac{p}{2}}r^{\alpha(n-p)+2}\sin(\phi) \ d\theta d\phi dr \notag \\
&=k^{p}\int_{1}^{4}\int_{0}^{\frac{\pi}{2}}\int_{0}^{2\pi}\Bigg[\frac{\cos^{2}(k\theta)}{k^{2}}\left(\frac{h'(r)^{2}\sin^{2}(\phi)}{(1+\alpha)^{2}}+(\frac{h(r)}{r})^{2}\cos^{2}(\phi)\right)\notag \\
&\qquad \qquad \qquad \qquad+\left(\frac{h(r)}{r}\right)^{2}\sin^{2}(k\theta)\Bigg]^{\frac{p}{2}}r^{\alpha(n-p)+2}\sin(\phi)\ d\theta d\phi dr \notag \\
&=k^{p}\int_{1}^{4}\int_{0}^{\frac{\pi}{2}}\int_{0}^{2\pi}\Bigg[\frac{\cos^{2}(u)}{k^{2}}\left(\frac{h'(r)^{2}\sin^{2}(\phi)}{(1+\alpha)^{2}}+\left(\frac{h(r)}{r}\right)^{2}\cos^{2}(\phi)\right)\notag\\
&\qquad \qquad \qquad \qquad+\left(\frac{h(r)}{r}\right)^{2}\sin^{2}(u)\Bigg]^{\frac{p}{2}}r^{\alpha(n-p)+2}\sin(\phi) \ du d\theta dr\notag
\end{align}
Therefore we have
$$\int_{\mathbb{R}^{3}}|\tilde{A}\nabla f_{k}|^{p}|x|^{\alpha(n-p)}\ dx=k^{p}\Bigg[ \int_{1}^{4}\int_{0}^{\frac{\pi}{2}}\int_{0}^{2\pi}\Bigg[\frac{h(r)}{r}\sin(u)\Bigg]^{p}r^{\alpha(n-p)+2}\sin(\phi)\ du d\phi dr+o(1) \Bigg]$$
A similar computation yields  $$\int_{\mathbb{R}^{3}}|\nabla f_{k}|^{p}|x|^{\alpha(n-p)}\ dx=k^{p}\Bigg[\int_{1}^{4}\int_{0}^{\frac{\pi}{2}}\int_{0}^{2\pi}\Bigg[\frac{h(r)}{r}\sin(u)\Bigg]^{p}r^{\alpha(n-p)+2}\sin(\phi)\ du d\phi dr +o(1) \Bigg]$$
Therefore
$$F(f_{k})=1+o(1), \text{ as } k\to \infty.$$
Combining this last estimate with Lemma 2.2, we get the conclusion of the Lemma.
\hfill$\Box$

Notice that with this Lemma, we have the proof of Theorem 1.2.
\begin{remark}
 One can see that this sequence of function $f_{k}$ always satisfies $$F(f_{k})\to 1 \text{ as } k\to \infty$$
for all $\alpha>-1$, but if  $\alpha<0$, we have that $A_{\alpha}=(1+\alpha)^{-t}>1$ thus the sequence $f_{k}$
in the case $\alpha<0$ is not optimizing and as we saw earlier, the optimizer is radially symmetric.
The sequence $f_{k}$ gains importance in the case $\alpha>0$ since $(\alpha+1)^{-t}<1$ hence there is a symmetry breaking
phenomena and the radially symmetric functions cannot be optimizers anymore.\\
On the other hand, by the Riemann-Lebesgue lemma,
$f_{k}\rightharpoonup 0$ as $k\to \infty$, hence we do not obtain an
optimizer in this case.
\end{remark}

\end{document}